\documentclass[12pt,leqno]{amsart} 

\usepackage{enumerate}
\usepackage{float} 
\usepackage{hyperref}

\usepackage{color}
\newlength{\fixboxwidth}
\newcommand{\Le}{{\mathcal L}}
\newcommand{\M}{{\mathcal M}}
\newcommand{\A}{{\mathcal A}}
\newcommand{\R}{{\mathbb R}}
\newcommand{\N}{{\mathbb N}}

\newcommand{\cost}{\operatorname{cost}}
\renewcommand{\rho}{{\varrho}}
\def\min{{\rm min}}

\def\a{{\alpha }} 
 
\def\e{{\varepsilon}}

\def\d{{\delta}} 
\def\phi{{\varphi}}  

\newcommand{\vol}{\operatorname{vol}}
\renewcommand{\rho}{{\varrho}}
\newcommand{\rad}{\mathcal R^{\alpha}}
\newcommand{\fad}{\mathcal F^\alpha}    
\newcommand{\ball}{{B^d}}

\newcommand{\wt}{\widetilde }

\newcommand{\mur}{\mu_\rho}

\newcommand{\abs}[1]{\left\vert #1 \right\vert} 
\newcommand{\norm}[2]{\left\Vert #1  \right\Vert _{#2}} 
 
\newcommand{\set}[1]{\left\{#1\right\}}
\newcommand{\expect}{\mathbf E}

\newcommand{\scalar}[2]{\left\langle #1,#2\right\rangle}

\theoremstyle{plain}
\newtheorem{theorem}{Theorem}
\newtheorem{lemma}[theorem]{Lemma}

\newtheorem{coro}[theorem]{Corollary}

\theoremstyle{definition} 
\newtheorem{rem}{Remark}

\floatstyle{boxed} \newfloat{alg}{htb}{}
\floatname{alg}{Algorithm}

\begin{document}

\title[Explicit error bounds for Markov Chain Monte Carlo]
{Explicit error bounds for lazy reversible\\ Markov Chain Monte Carlo} 

\author{Daniel Rudolf} 
\address{ Friedrich Schiller University Jena, 
Mathem. Institute,
Ernst-Abbe-Platz 2, 
D-07743 Jena, Germany}
\email{ruda@minet.uni-jena.de}
\date{Version: \today}

\begin{abstract}
We prove explicit, i.e., non-asymptotic, error bounds 
for Markov Chain Monte Carlo methods,
such as the Metropolis algorithm.
The problem is to compute the expectation (or integral)
of $f$ with respect to a measure $\pi$ which can be given
by a density $\rho$ with respect to another measure.
A straight simulation of the desired distribution by a random number generator
 is in general not possible.
Thus it is reasonable to use Markov chain sampling with a burn-in. 
We study such an algorithm and 
extend the analysis of Lovasz and Simonovits (1993)
to obtain an explicit error bound.
\end{abstract}

\maketitle

\section{Problem description, Introduction} \label{problem_descr} 

The paper deals with numerical integration based on Markov chains. 
The main goal is to approximate an integral of the following form
\begin{equation}  \label{int}
S(f):=\int_\Omega f(x) \, \pi(dx),  
\end{equation}
where $\Omega$ is a given set and $\pi$ a probability measure. 
In addition we assume that an oracle which computes function values of $f$ is provided.
We generate a Markov chain $X_1,X_2,\dots$ with transition kernel $K$, having $\pi$ as its stationary distribution.
After a certain burn-in time there is an average computation over the generated sample (Markov chain steps). 
For a given function $f$ and burn-in time, say $n_0$, we get as approximation
\begin{equation*}  \label{sum}
	S_{n,n_0}(f):=\frac{1}{n} \sum_{j=1}^n f(X_{j+n_0}).
\end{equation*} 

This Markov chain Monte Carlo method (MCMC) for approximating the expectation 
plays a crucial role
in numerous applications, especially in statistical physics, in statistics,
and in financial mathematics.
Certain asymptotic error bounds are known, which can be proved via
isoperimetric inequalities, the Cheeger inequality and estimates of eigenvalues, see \cite{sokal,mathe1,mathe2}.
Here in contrast, we determine an explicit error bound for $S_{n,n_0}$.
The individual error of such a method $S_{n,n_0}$ and a function $f$ is measured in mean square sense, i.e.,
\begin{equation*}
	\label{error_ind}
  e(S_{n,n_0}	,f):=\left( \expect \abs{	S_{n,n_0}(f)-S(f)}^2 \right)^{1/2}.	
\end{equation*}

Now an outline of the structure of the paper 
and the main results
is given.
Section~\ref{basics} contains
the used notation and repeats some relevant statements.
An introduction of the idea of laziness is given in Section~\ref{lazy}, 
where also the conductance concept and a convergence property of the chain is presented.
It is useful for getting results to restrict ourself to Markov chains which have a positive conductance 
and where the initial distribution $\nu$, for obtaining the first time step, has a bounded density with respect to $\pi$.
Section~\ref{error_bounds} 
contains the new results.     
Let $\phi$ be the conductance of the
underlying chain. After a burn-in 
\[ 
n_0 \geq \frac{\log\left(\norm{\frac{d\nu}{d\pi}}{\infty}\right)}{\phi^2}
\quad 
\mbox{the error obeys} 
\quad
e(S_{n,n_0},f)\leq \frac{10}{\phi\cdot\sqrt{n}}\norm{f}{\infty}.
\] 
This implies immediately that the number $n+n_0$ of time steps which are needed for an error $\e$, 
can be bounded by 
\[
\left\lceil \frac{\log\left(\norm{\frac{d\nu}{d\pi}}{\infty}\right)}{\phi^2} \right\rceil 
+    \left\lceil \frac{100\norm{f}{\infty}^2}{\phi^2\cdot\e^2}\right\rceil.
\]
All results are in a general framework,
such that after an adaption it is possible to apply the theory in different 
settings e.g. discrete state space or continuous one.   
In Section~\ref{appl} we pick up a problem considered in \cite{novak}.
There the authors use the Metropolis algorithm for approximating an integral over the $d$ dimensional unit ball $\ball\subset \R^d$
with respect to an unnormalised density.
The strict positive density is notated by $\rho$ and moreover 
we assume that it is logconcave and $\a$ is the Lipschitz constant of $\log\rho$.
Let $\delta >0$ and $B(x,\d)$ be the ball with radius $\d$ around $x$. 
Then we suggest the method described in Algorithm~\ref{mcmc} (see page \pageref{mcmc}) for the approximation of 
\[
S(f)=S(f,\rho)=\frac{\int_\ball f(x)\rho(x) dx}{\int_\ball \rho(x) dx}.
\]

\noindent
It is shown that for $\d=\min\set{1/\sqrt{d+1},1/\a}$ the error obeys
\[
e(S_{n,n_0}^\d,f)\leq 
8000 \frac{\sqrt{d+1} \max\set{ \sqrt{d+1},\a}}{\sqrt{n}}\norm{f}{\infty},
\] 
where the burn-in time $n_0$ is chosen larger than $1280000\cdot\a (d+1) \max\set{d+1,\a^2}$. 
\begin{alg}[ht]
\begin{flushleft}
\textbf{\quad Algorithm: $S_{n,n_0}^\d(f,\rho)$} 
\end{flushleft}
\hrule
\begin{enumerate}
  \item 	\texttt{choose $X_1$ randomly on $\ball$;}\\
  \item   \texttt{for $i=1,\dots,n+n_0$ do\\
  \begin{itemize}
 						 	\item if $rand()< 1/2$ then $X_{i+1}:= X_i$;\\
 						 	\item else 
 						 	\begin{itemize}
 						 	  \item[-] choose $Y\in B(X_i,\delta)$ uniformly;\\
  							\item[-] if $Y \notin \ball$ then $X_{i+1}:= X_i$;\\
  							\item[-] if $Y\in \ball$ and $\rho(Y) \ge \rho(X_i)$ then
                 $X_{i+1}:=Y$;\\			
  							\item[-] if $Y\in \ball$ and $\rho(Y) \le \rho(X_i)$ then\\ 
  							\item[-] \quad $X_{i+1}:=Y$ with Prob $\rho(Y)/\rho(X_i)$ and \\ 
  							\item[-]	\quad$X_{i+1}:=X_i$  with Prob $1-\rho(Y)/\rho(X_i)$.
  						\end{itemize}
  \end{itemize}}
  \item   \texttt{Return:     \\  
					\[
					S_{n,n_0}^\d(f,\rho):=\frac{1}{n} \sum_{j=1}^n f(X_{j+n_0}).
					\]}
\end{enumerate}
\caption{Metropolis algorithm for $S(f,\rho)$}
\label{mcmc}
\end{alg}

It is worth pointing out that the number of time steps which we use for sampling
behaves polynomial in the dimension and also polynomial
in the Lipschitz constant $\a$ of the densities.
As already mentioned the same integration problem was studied in \cite{novak}.
The authors asked whether the problem is tractable.
That means the number of function evaluation 
to obtain an error smaller than $\e$  can be polynomially bounded 
by the dimension and the Lipschitz constant. 
So we give a positive answer; the problem is tractable, at least if we consider bounded integrands $f$.

\goodbreak

\section{Notation and basics} \label{basics}
In this section we explain the most important facts and definitions which we are going to use in the analysis.
For introductory literature to general Markov chains we refer the reader to \cite{tweed}, \cite{Nummelin} or \cite{revuz}.
Throughout this study we assume that $(\Omega,\A)$ is a measurable countably generated space.
Then we call $K: \Omega \times \A \to [0,1]$ Markov kernel or transition kernel if   
\begin{enumerate}[(i)]
 \item for each $x\in\Omega$ the mapping $A\in\A \mapsto K(x,A)$ induces a probability measure on $\Omega$, 
 \item for each $A\in\A$ the mapping $x\in\Omega \mapsto K(x,A)$ is an $\A$-measurable real function.
\end{enumerate}

\goodbreak

In addition $\M=(\Omega,\A,\{ K(x,\cdot): x\in\Omega \})$ is the associated Markov scheme. This notation is taken from \cite{lova_simo1}.
A Markov chain $X_1,X_2,\dots$ is given through a Markov scheme $\M$ and a start distribution $\nu$ on $\Omega$. The transition kernel $K(x,A)$ of the Markov chain describes the probability of getting from $x\in\Omega$ to $A\in\A$ in one step. 
Another important assumption is that the given distribution $\pi$ is stationary concerning the considered Markov chain, i.e. for all $A\in\A$
\begin{equation*}
	\pi(A)=\int_\Omega K(x,A) \pi(dx).
\end{equation*}  
Roughly spoken that means: Choosing the starting point with distribution $\pi$, then after one step we have the same distribution as before. Another similar but stronger restriction of the chain is reversibility. A Markov scheme is reversible with respect to $\pi$ if for all $A, B\in\A$
\begin{equation*}
 \int_B K(x,A) \pi(dx)=\int_A K(x,B) \pi(dx).	
\end{equation*} 
The next outcome is taken from \cite{lova_simo1}. 
But it is not proven there so we will give an idea of the proof.
\begin{lemma}  \label{f_eq}
Let $\M$ be a reversible Markov scheme and let $F:\Omega \times \Omega \to \R$ be integrable. Then 
\begin{equation}   \label{F_x_y}
\int_\Omega \int_\Omega F(x,y)\;K(x,dy)\pi(dx)=\int_\Omega \int_\Omega F(y,x)\;K(x,dy)\pi(dx).
\end{equation} 

\end{lemma}
\begin{proof}
The result is shown using a standard technique of integration theory. 
Since the Markov scheme is reversible we have
\[
\int_\Omega \int_\Omega I_{A\times B} (x,y) K(x,dy) \pi(dx)=\int_\Omega \int_\Omega I_{A\times B} (y,x) K(x,dy) \pi(dx)
\]
for $A,B\in\A$. Having finished this we develop the equality of the integrals for an arbitrary set $C\in \A \otimes \A$, 
where $\A \otimes \A$ is the product $\sigma$-algebra of $\A$ with itself. 
This is an application of the Dynkin system theorem. 
Then we consider the case where $f$ is a simple function, which is straightforward. 
The next step is to obtain the equality for positive function and after that extending the result to general integrable ones. 
\end{proof}

\begin{rem} \label{rem_f}
If we have a Markov scheme, which is not necessarily reversible but has a stationary distribution the following holds true
\[
S(f)=\int_\Omega f(x) \pi(dx)= \int_\Omega\int_\Omega f(y) K(x,dy) \pi(dx) ,
\]
where $f:\Omega\to\R$ is integrable. 
This can be seen easily by using the same steps as in the proof of Lemma~\ref{f_eq}.
\end{rem}

By $K^n(x,\cdot)$ we denote the $n$-step transition probabilities and we have for $x\in\Omega$, $A\in\A$ that
\[
K^n(x,A):=\int_\Omega K^{n-1}(y,A) K(x,dy)=\int_{\Omega}K(y,A)K^{n-1}(x,dy).
\]
This again constitutes a transition kernel of a Markov chain
sharing the invariant distribution and reversibility with the original one.
Thus the outcomes of Lemma~\ref{f_eq} and Remark~\ref{rem_f} also hold for
the $n$-step transition probabilities i.e.
\begin{equation}  \label{F_x_y_k}
\int_\Omega \int_\Omega F(x,y)\,K^n(x,dy)\,\pi(dx)=\int_\Omega \int_\Omega F(y,x)\;K^n(x,dy)\,\pi(dx).
\end{equation}


Now we define for a Markov scheme $\M$ a nonnegative 
operator $P: L_\infty(\Omega,\pi) \to L_\infty(\Omega,\pi)$ by
\[
(Pf)(x)=\int_\Omega f(y) K(x,dy).
\]
(Nonnegative means: if $f\geq0$ then $Pf\geq0$.)
This operator is called Markov or transition operator concerning a Markov scheme $\M$ and describes the expected value of $f$ after one step with the Markov chain from $x\in\Omega$.
The expected value of $f$ from $x\in\Omega$ after $n$-steps with the Markov chain is given as 
\[
(P^nf)(x)= \int_\Omega f(y) K^n(x,dy).
\] 
Let us now consider $P$ on the Hilbert space $L_2(\Omega,\pi)$ 
and $\scalar{f}{g}=\int_\Omega f(x)g(x)\,\pi(dx)$ denotes the canonical scalar product. 
Notice that the considered function space is chosen according to the invariant measure.
Then we have with Lemma~\ref{f_eq}
\begin{equation} \label{greater_zero}
\scalar{f}{f}\pm\scalar{f}{Pf}=\frac{1}{2}\int_\Omega\int_\Omega (f(x)\pm f(y))^2 K(x,dy) \pi(dx)  \geq0.
\end{equation}
From a functional analysis point of view that means $\norm{P}{L_2\to L_2}\leq1$.
It is straightforward to show that $\norm{P^n}{L_p\to L_p}\leq1$ for $p=1,2$ or $\infty$ and $n\in\N$.


Let $X_1,X_2,\dots$ be the result of a reversible Markov chain.
The expectation of the chain with starting distribution $\nu=\pi$ and Markov kernel $K$
from scheme $\M$ is denoted by $\expect_{\pi,K}$.
Then we get for $f\in L_2(\Omega,\pi)$
\begin{align}  \label{expect_scalar}
\notag & \expect_{\pi,K}(f(X_i)) =\expect_{\pi,K}(f(X_0))=\scalar{1}{f}=S(f),\\
\notag & \expect_{\pi,K}(f(X_i)^2) =\expect_{\pi,K}(f(X_0)^2)=\scalar{f}{f}=S(f^2),\\
& \expect_{\pi,K}(f(X_i)f(X_{j})) =\expect_{\pi,K}(f(X_0)f(X_{\abs{i-j}}))=\scalar{f}{P^{\abs{i-j}} f}.
\end{align}

The assumption that the initial distribution is the stationary one makes the calculation easy. 
In the general case, where the starting point is chosen by a given probability distribution $\nu$, 
we obtain for $i\leq j$ and functions $f\in L_2(\Omega,\pi)$
\begin{align*}
& \expect_{\nu,K}(f(X_i)) =\int_\Omega P^i f(x) \nu(dx),\\
& \expect_{\nu,K}(f(X_i)f(X_{j})) = \int_\Omega P^i(f(x)P^{j-i}f(x)) \nu(dx).
\end{align*}
It is easy to verify with \eqref{F_x_y} that $P$ is self-adjoint as acting on $L_2(\Omega,\pi)$. 
In the next part we are going to get one more convenient characteristic of $P$ 
under some additional restrictions.

\section{Laziness and Conductance}  \label{lazy}
An introduction to laziness and a more detailed view on the conductance 
is given in \cite{lova_simo1}. Most results which we are going to mention here
are taken from this reference.  
A Markov scheme $\M=(\Omega,\A,\{ K(x,\cdot): x\in\Omega \})$ is called lazy if 
$K(x,\set{x})\geq1/2$ for all $x\in\Omega$. 
This means the chain stays at least with probability $1/2$ in the current state.
Notice that the resulting chain from Algorithm~\ref{mcmc} (see page \pageref{mcmc}) is lazy because of line three. 
The crucial fact for slowing down is to deduce that the associated Markov operator $P$ is positive semidefinite. 
Therefore we study only lazy chains.
This is formalized in the next Lemma.

\begin{lemma}  \label{pos_def}
Let $\M$ be a lazy, reversible Markov scheme then we have
for $f\in L_2(\Omega,\pi)$ 
\begin{equation}  \label{pos_def_eq}
\scalar{Pf}{f}\geq0.
\end{equation}  
\end{lemma}

\begin{proof}
We consider another Markov scheme $\wt{\M}:=(\Omega,\A,\{\wt{K}(x,\cdot):x\in\Omega\})$, where 
$\wt{K}(x,A)=2K(x,A)-I(x,A)$ with 
\[
I(x,A)=\begin{cases} 1 & x\in A\\
                     0 & x\in A^c
       \end{cases}
\]
for all $A\in\A$. To verify, that $\wt{K}$ is again a transition kernel we need $K(x,\set{x})\geq1/2$. 
The reversibility condition for $\wt{\M}$ holds, since scheme $\M$ is reversible.
The Markov operator of $\wt{\M}$ is given by $\wt{P}=(2P-I)$, where $I$ is the identity.
Since we established reversibility of the new scheme we obtain by applying Lemma~\ref{f_eq} equality \eqref{greater_zero} for $\wt{P}$.
So it is true that
\[
-\scalar{f}{f}\leq \scalar{(2P-I)f}{f} \leq \scalar{f}{f}.
\]
Now let us consider 
\[
\scalar{Pf}{f}=\frac{1}{2}\scalar{f}{f}+\frac{1}{2}\scalar{(2P-I)f}{f}\geq0,
\]
such that the claim is proven.
\end{proof}


Having finished this, we can turn to the conductance of the Markov chain.
For a Markov scheme $\M=(\Omega,\A,\{ K(x,\cdot): x\in\Omega \})$, 
which is not necessarily lazy, it is defined by
\[
  \phi(K,\pi)=\inf_{0<\pi(A)\leq1/2}\frac{\int_A K(x,A^c)\pi(dx)}{\pi(A)},
\]
where $\pi$ is a stationary distribution.
The numerator of the conductance describes the 
probability of leaving $A$ in one step, 
where the starting point is chosen by $\pi$.
An important requirement for the following is 
that the scheme has a positive conductance, 
since the next result is not useful otherwise.

\begin{lemma} \label{mixing_lem}
Let $\M$ be a lazy, reversible Markov scheme and let $\nu$ be the initial distribution. 
Furthermore we assume that the probability distribution $\nu$ has a bounded density function $\frac{d\nu}{d\pi}$ with respect to $\pi$. 
Then for $A\in\A$
we obtain
\begin{equation}  \label{mixing_res}
\abs{\int_\Omega K^j(x,A)\,\nu(dx)-\pi(A)}
\leq \sqrt{\norm{\frac{d\nu}{d\pi}}{\infty}} \left( 1-\frac{\phi(K,\pi)^2}{2}  \right)^j.
\end{equation} 
\end{lemma}
\begin{proof}
Look at the result of \cite[Corollary~1.5, p. 372]{lova_simo1} and translate it in our notation. 
\end{proof}

\begin{rem}
The left hand side of \eqref{mixing_res} can be transformed as follows
\begin{align*}
& \int_\Omega K^j(x,A)\:\nu(dx)-\pi(A) = \int_\Omega \int_A K^j(x,dy)\: \frac{d\nu}{d\pi}(x)\: \pi(dx)-\pi(A)\\
& \underset{\eqref{F_x_y_k}}{=}  \int_A \int_\Omega \frac{d\nu}{d\pi}(y)\: K^j(x,dy)\:  \pi(dx)
														-\int_A \int_\Omega \frac{d\nu}{d\pi}(y) \:\pi(dy)\: \pi(dx)\\
& = \int_A \int_\Omega \frac{d\nu}{d\pi}(y)\: (K^j(x,dy)-\pi(dy))\:  \pi(dx).
\end{align*}
Now it is clear that with Lemma~\ref{mixing_lem} for $A\in\A$
\begin{equation}  \label{mixing_res2}
\abs{\int_A \int_\Omega \frac{d\nu}{d\pi}(y)\: (K^j(x,dy)-\pi(dy))\:  \pi(dx)}
\leq \sqrt{\norm{\frac{d\nu}{d\pi}}{\infty}} \left( 1-\frac{\phi(K,\pi)^2}{2}  \right)^j.
\end{equation}
\end{rem}
\begin{rem}
Observe, that we got a bound for the speed of convergence to stationarity of the considered Markov chain. 
Once more it is possible to estimate the right hand side of \eqref{mixing_res2}, in detail 
\begin{equation}  \label{exp_j}
\sqrt{\norm{\frac{d\nu}{d\pi}}{\infty}} \left( 1-\frac{\phi(K,\pi)^2}{2}  \right)^j
\leq\sqrt{\norm{\frac{d\nu}{d\pi}}{\infty}} \exp\left[-j\frac{\phi(K,\pi)^2}{2}\right]
\end{equation}
holds true.
\end{rem}

To use the conductance we need a connection to the operator $P$. 
This is given in form of the so called Cheeger inequality. 
Before we are going to state this conclusion in a slightly different formulation 
we define a subset of $L_2(\Omega,\pi)$ as follows
\[
L_2^0=L_2^0(\Omega,\pi):=\set{f\in L_2(\Omega,\pi): S(f)=0}.
\] 
\goodbreak
\begin{lemma}[Cheeger's inequality] \label{cheeger_lem}
Let $\M$ be a reversible Markov scheme with conductance $\phi(K,\pi)$.
Then for $g\in L_2^0$
\begin{equation} \label{cheeger_eq}
 \scalar{P^jg}{g}\leq \left( 1-\frac{\phi(K,\pi)^2}{2} \right)^j \norm{g}{2}^2.
\end{equation}
\end{lemma}
\begin{proof}
See \cite[Corollary~1.8, p. 375]{lova_simo1}.
\end{proof}
\begin{rem}
There are many other references where the convergence rate of Markov chains to stationarity is studied, see e.g.  \cite{jerrum,dia_stroock_eigen,rosenthal_minor,rosenthal}.
One approach is to bound the second eigenvalue of the operator $P$. 
The relation between the eigenvalue and the conductance of a Markov chain is given by Cheeger's inequality (see Lemma~\ref{cheeger_lem}).
In this context the laziness condition shifts the spectrum of the Markov operator $P$ restricted to $L_2^0$ from $(-1,1)$ by the transformation described in Lemma~\ref{pos_def}
to $(0,1)$ i.e. the second eigenvalue is always positive.
\end{rem}

\section{Error bounds} \label{error_bounds}
This section contains the main result and its proof.
At first we are going to repeat an already known finding, which is used to show an explicit error bound for 
a general Markov scheme with initial probability distribution $\nu$. Most arguments to obtain that result are from \cite{lova_simo1} and \cite{mathe1}. \\

The next conclusion considers an algorithm 
under the assumption that the starting point is chosen according to the stationary distribution.
So a preliminary burn-in period is not necessary anymore since we are already 
at the invariant distribution.

\begin{theorem}  \label{lovasz_int}
Let $\M$ be a lazy, reversible Markov scheme with stationary distribution $\pi$, let $X_1,X_2,\dots$ be a Markov chain generated by $\M$ with initial distribution $\pi$. Let $f\in L_2(\Omega,\pi)$, $S(f)=\int_\Omega f(x) \pi(dx)$ and $S_n(f):=S_{n,0}(f)=\frac{1}{n}\sum_{j=1}^n f(X_{j})$. Then  we obtain
\[
e(S_n,f)^2=\expect_{\pi,K}\abs{S(f)-S_n(f)}^2 \leq \frac{4}{\phi(K,\pi)^2\cdot n}  \norm{f}{2}^2.
\]
\end{theorem}
\begin{rem}
This proof is again taken from \cite[Theorem~1.9, p. 375]{lova_simo1}. 
Since it is very important in our analysis and 
because of the slightly different notation we will repeat it.
\end{rem}
\begin{proof}
 Let $g:=f-S(f)$, such that $g\in L_2^0$. 
 Then we have with  Lemma~\ref{pos_def}, Lemma~\ref{cheeger_lem} and $\norm{g}{2}\leq\norm{f}{2}$ that
 \begin{align*}
  &\expect_{\pi,K}\abs{S(f)-S_n(f)}^2=\expect_{\pi,K} \abs{\frac{1}{n}\sum_{j=1}^n g(X_{j})}^2\\
  &=\frac{1}{n^2} \sum_{j=1}^n \sum_{i=1}^n \expect_{\pi,K}(g(X_{j})g(X_{i}))
   \underset{\eqref{expect_scalar}}{=} \frac{1}{n^2} \sum_{j=1}^n \sum_{i=1}^n \expect_{\pi,K}(g(X_0)g(X_{\abs{i-j}}))\\
  &=\frac{1}{n^2}\left(n\scalar{g}{g} + \sum_{k=1}^{n-1}2(n-k)\scalar{P^kg}{g} \right) \\
  &\leq\frac{1}{n^2} \sum_{k=0}^{n-1} 2(n-k)\scalar{P^kg}{g} 
   \underset{\eqref{pos_def_eq}}{\leq} \frac{2}{n} \sum_{k=0}^\infty \scalar{P^kg}{g}
   \underset{\eqref{cheeger_eq}}{\leq}  \frac{2}{n} \sum_{k=0}^\infty \left (1-\frac{\phi(K,\pi)^2}{2} \right)^k \norm{g}{2}^2\\
  &=\frac{4}{\phi(K,\pi)^2\cdot n} \norm{g}{2}^2
   \leq \frac{4}{\phi(K,\pi)^2\cdot n}  \norm{f}{2}^2.
 \end{align*}
Notice that laziness is essentially used by applying $\scalar{P^kg}{g}\geq0$ in the second inequality.
\end{proof}
Let us consider the more general case, where the initial distribution is not the stationary one. 
In the next statement a relation between the error of starting with $\pi$ and the error
of starting not with the invariant distribution is established.
\begin{lemma} \label{connect_lem}
Let $\M$ be a reversible Markov scheme with stationary distribution $\pi$, let $X_1,X_2,\dots$ 
be a Markov chain generated by $\M$ with initial distribution $\nu$. 
Let $\frac{d\nu}{d\pi}$ be a bounded density of $\nu$ with respect to $\pi$. Then we get for $g:=f-S(f)\in L_2^0$
\begin{multline}  \label{connection}
\expect_{\nu,K} \abs{S(f)-S_{n,n_0}(f)}^2 = \expect_{\pi,K} \abs{S(f)-S_n(f)}^2\\ 
+ \frac{1}{n^2}\sum_{j=1}^{n} \int_\Omega \int_\Omega \frac{d\nu}{d\pi}(y) \left(K^{n_0+j}(x,dy)-\pi(dy)\right)\;g(x)^2 \; \pi(dx)\\
+ \frac{2}{n^2} \sum_{j=1}^{n-1} \sum_{k=j+1}^n
  \int_\Omega \int_\Omega \frac{d\nu}{d\pi}(y) \;\left(K^{n_0+j}(x,dy)-\pi(dy)\right)\;g(x) P^{k-j}g(x)\;\pi(dx).
\end{multline}
\end{lemma}
\begin{proof}
It is easy to see, that
\begin{align*}
& \expect_{\nu,K} \abs{S(f)-S_{n,n_0}(f)}^2=\frac{1}{n^2} \sum_{j=1}^n \sum_{i=1}^n \expect_{\nu,K} (g(X_{n_0+j})g(X_{n_0+i}))\\
&= \frac{1}{n^2} \sum_{j=1}^n \int_\Omega P^{n_0+j}g(x)^2\; \nu(dx)+ \frac{2}{n^2} \sum_{j=1}^{n-1} \sum_{k=j+1}^n \int_\Omega P^{n_0+j}(g(x) P^{k-j}g(x))\;\nu(dx).
\end{align*}
For every function $h\in L_2(\Omega,\pi)$ and $i\in\N$ under applying \eqref{F_x_y_k} the following transformation holds true
\begin{align*} 
& \int_\Omega P^i h(x)\, \nu(dx) 
	= \int_\Omega \int_\Omega h(y)\, K^i(x,dy)\, \frac{d\nu}{d\pi}(x)\, \pi(dx)\\
& \underset{\eqref{F_x_y_k}}{=} 
	\int_\Omega \int_\Omega \frac{d\nu}{d\pi}(y)  K^i(x,dy)\,h(x) \, \pi(dx)\\
& = \int_\Omega h(x)\, \pi(dx) 
	+ \int_\Omega \int_\Omega \frac{d\nu}{d\pi}(y)  \left(K^i(x,dy)-\pi(dy)\right) h(x) \, \pi(dx)\\
& \underset{\eqref{F_x_y_k}}{=} 
    \int_\Omega P^i h(x)\pi(dx) 
	+ \int_\Omega \int_\Omega \frac{d\nu}{d\pi}(y)  \left(K^i(x,dy)-\pi(dy)\right) h(x) \, \pi(dx).
\end{align*}
Using this in the above setting formula \eqref{connection} is shown.
\end{proof}
The next finding is also a helpful tool to prove the main result of this paper. 
It modifies the convergence property, which is described in Lemma~\ref{mixing_lem}, 
such that we are able to use it in the considered context.
\begin{lemma} \label{Omega_vorz}
Let $\M$ be a lazy, reversible Markov scheme with stationary distribution $\pi$, let $\nu$ be the initial distribution with bounded density $\frac{d\nu}{d\pi}$ of the related Markov chain. Then we obtain for $h\in L_\infty(\Omega,\pi)$ and $j\in\N$ 
\[
\abs{\int_\Omega \int_\Omega \frac{d\nu}{d\pi}(y)  \left(K^j(x,dy)-\pi(dy)\right)h(x) \, \pi(dx)} 
\leq 4 \norm{h}{\infty} \sqrt{\norm{\frac{d\nu}{d\pi}}{\infty}}\left( 1-\frac{\phi(K,\pi)^2}{2} \right)^j.
\]
\end{lemma} 
\begin{proof}
At first we define $p_j(x):=\int_\Omega \frac{d\nu}{d\pi}(y)  \left(K^j(x,dy)-\pi(dy)\right)$.
With the standard proof technique of integration theory it is easy to see that the measurability of 
the density and the kernel can be carried over to $p_j$. 
Now we consider the positive and negative parts of the functions $h$ and $p_j$. To formalize this we use
\begin{align*}
\Omega^{+}_+& :=\set{x\in\Omega: p_j(x)\geq0,\;h(x)\geq0},\\
\Omega^{+}_-& :=\set{x\in\Omega: p_j(x)\geq0,\;h(x)<0},\\
\Omega^{-}_+& :=\set{x\in\Omega: p_j(x)<0,\;h(x)\geq0},\\
\Omega^{-}_-& :=\set{x\in\Omega: p_j(x)<0,\;h(x)<0}.
\end{align*}
These subsets of $\Omega$ are all included in the $\sigma$-algebra $\A$, 
since $p_j$ and $h$ are measurable functions.
So applying \eqref{mixing_res2} leads to the following upper bound
\begin{align*}
 &\abs{\int_\Omega p_j(x)h(x)\,\pi(dx)} 
  \leq \abs{\int_{\Omega^+_+} p_j(x)h(x)\,\pi(dx)} + \abs{\int_{\Omega^+_-} p_j(x)h(x)\,\pi(dx)} \\
& \qquad +  \abs{\int_{\Omega^{-}_+}p_j(x)h(x)\,\pi(dx)} +  \abs{\int_{\Omega^{-}_-}p_j(x)h(x)\,\pi(dx)}\\
&   \leq \norm{h}{\infty} \abs{\int_{\Omega^{+}_+}p_j(x)\,\pi(dx)} 
  + \norm{h}{\infty} \abs{\int_{\Omega^{+}_-}p_j(x)\,\pi(dx)} \\
&	\qquad+ \norm{h}{\infty}  \abs{\int_{\Omega^{-}_+}p_j(x)\,\pi(dx)} 
  + \norm{h}{\infty} \abs{\int_{\Omega^{-}_-}p_j(x)\,\pi(dx)}\\
&\underset{\eqref{mixing_res2}}{\leq} 4 \norm{h}{\infty} \sqrt{\norm{\frac{d\nu}{d\pi}}{\infty}}\left( 1-\frac{\phi(K,\pi)^2}{2} \right)^j.
\end{align*}
\end{proof}
Now all results are available to obtain 
our main error bound for the MCMC method $S_{n,n_0}$.
\begin{theorem}  \label{error_thm}
Let $X_1,X_2,\dots$ be a lazy, reversible Markov chain, defined by the scheme $\M$ and the initial distribution $\nu$.
Let the initial distribution have a bounded density $\frac{d\nu}{d\pi}$ with respect to $\pi$. Let 
$S_{n,n_0}(f)=\frac{1}{n}\sum_{j=1}^n f(X_{n_0+j})$ be the approximation of $S(f)=\int_\Omega f(x) \pi(dx)$, where $f\in L_\infty(\Omega,\pi)$. Then
\begin{equation*}  \label{error}
e(S_{n,n_0},f)\leq \frac{2\sqrt{1+24\sqrt{\norm{\frac{d\nu}{d\pi}}{\infty}}
									\exp\left[-n_0\frac{\phi(K,\pi)^2}{2}\right]}}{\phi(K,\pi)\cdot \sqrt{n}}\norm{f}{\infty}.
\end{equation*}
\end{theorem}

\begin{proof}
By Lemma~\ref{connect_lem} and Lemma~\ref{Omega_vorz} where $g:=f-S(f)$ we have
\begin{align*}  
\expect_{\nu,K} &\abs{S(f)-S_{n,n_0}(f)}^2 \leq \expect_{\pi,K} \abs{S(f)-S_n(f)}^2\\ 
&\qquad\qquad+ \frac{4\norm{g}{\infty}^2}{n^2}\sum_{j=1}^{n} 
      \sqrt{\norm{\frac{d\nu}{d\pi}}{\infty}}\left( 1-\frac{\phi(K,\pi)^2}{2} \right)^{j+n_0}\\
&\qquad\qquad+ \frac{8\norm{g}{\infty}^2}{n^2} \sum_{j=1}^{n-1} \sum_{k=j+1}^n \norm{P^{k-j}}{L_\infty\to L_\infty} 
  \sqrt{\norm{\frac{d\nu}{d\pi}}{\infty}}\left( 1-\frac{\phi(K,\pi)^2}{2} \right)^{j+n_0}.
\end{align*}
For an easier notation we define 
\begin{equation}  \label{eps_0}
\e_0:=\sqrt{\norm{\frac{d\nu}{d\pi}}{\infty}}\exp\left[-n_0\frac{\phi(K,\pi)^2}{2}\right].
\end{equation}
Taking \eqref{exp_j} and \eqref{eps_0} into account the following transformation is true 
\begin{align*}  
\expect_{\nu,K} \abs{S(f)-S_{n,n_0}(f)}^2 & \leq \expect_{\pi,K} \abs{S(f)-S_n(f)}^2 
+ \frac{4\,\e_0\norm{g}{\infty}^2}{n^2}\sum_{j=1}^{n} 
      \left( 1-\frac{\phi(K,\pi)^2}{2} \right)^{j}\\
&\quad+ \frac{8\,\e_0\norm{g}{\infty}^2}{n^2} \sum_{j=1}^{n-1} \sum_{k=j+1}^n \norm{P^{k-j} }{L_\infty\to L_\infty}
  \left( 1-\frac{\phi(K,\pi)^2}{2} \right)^{j}.
\end{align*}
With the geometric series and $\norm{P^i}{L_\infty\to L_\infty}\leq1$ for all $i\in\N$ we get
\begin{align*}  
\expect_{\nu,K} \abs{S(f)-S_{n,n_0}(f)}^2 & \leq \expect_{\pi,K} \abs{S(f)-S_n(f)}^2 
+ \frac{8\,\e_0\norm{g}{\infty}^2}{\phi(K,\pi)^2\cdot n^2}\\
&\quad\quad\quad+ \frac{8\,\e_0\norm{g}{\infty}^2}{n^2} \sum_{j=1}^{n-1} (n-j) \left( 1-\frac{\phi(K,\pi)^2}{2} \right)^{j}\\
&\leq \expect_{\pi,K} \abs{S(f)-S_n(f)}^2 
+ \frac{16\,\e_0\norm{g}{\infty}^2}{\phi(K,\pi)^2\cdot n}+\frac{8\,\e_0\norm{g}{\infty}^2}{\phi(K,\pi)^2\cdot n^2}\\
&\leq \expect_{\pi,K} \abs{S(f)-S_n(f)}^2 + \frac{24\,\e_0\norm{g}{\infty}^2}{\phi(K,\pi)^2\cdot n}
.
\end{align*}

After applying Theorem~\ref{lovasz_int} and using $\norm{f}{2}^2\leq\norm{f}{\infty}^2$, $\norm{g}{\infty}^2\leq4\norm{f}{\infty}^2$ everything is proven.
\end{proof}


The major difference between the new error bound of Theorem~\ref{error_thm} and the already 
known from Theorem~\ref{lovasz_int} is that the unrealistic assumption 
to sample from the stationary distribution $\pi$ for the first time step
is weakened. It came out that for a certain burn-in time $n_0$ a very similar 
upper bound holds true, if the initial distribution $\nu$ has a bounded density
with respect to $\pi$. 
A further estimation yields the next conclusion.
\begin{coro}  
Let $X_1,X_2,\dots$ be a lazy, reversible Markov chain. The initial distribution $\nu$ has a bounded density $\frac{d\nu}{d\pi}$ with respect to $\pi$. 
Then for $f\in L_\infty(\Omega,\pi)$ and $S_{n,n_0}(f)=\frac{1}{n}\sum_{j=1}^n f(X_{j+n_0})$ after a burn-in 
\begin{equation} \label{error_case}
n_0\geq \frac{\log\left( \norm{\frac{d\nu}{d\pi}}{\infty} \right)}{\phi(K,\pi)^2}
\quad \mbox{the error obeys} \quad
e(S_{n,n_0},f)\leq \frac{10}{\phi(K,\pi)\cdot\sqrt{n}}\norm{f}{\infty}.
\end{equation}
\end{coro}

If we denote by $\cost(f,\e)$ the number $n+n_0$ of time steps that are 
needed for an optimal algorithm to solve \eqref{int} within an error $\e$,
then we can also write
\begin{equation*} \label{cost}
\cost(f,\e)
\leq \left\lceil \frac{\log\left(\norm{\frac{d\nu}{d\pi}}{\infty}\right)}{\phi(K,\pi)^2} \right\rceil 
+    \left\lceil \frac{100\norm{f}{\infty}^2}{\phi(K,\pi)^2\cdot\e^2}\right\rceil.
\end{equation*}
Roughly spoken that means if we control the conductance
of the underlying Markov chain, then we also control 
the error. So we should look for lower bounds 
of the conductance to obtain upper estimations of the error.

The next task is to apply the received results for an 
explicit example where we can use \eqref{error_case}.

\goodbreak

\section{Application} \label{appl}
For working with the above presented theory we need a lazy and reversible Markov chain.
In the following a construction for a reversible Markov scheme, using the Metropolis algorithm, 
is provided. After having this scheme we make it lazy and carry 
the conductance properties over to the new chain. 
This laziness is easily obtained by pasting a coin tossing step, 
where we accept the new state when head occurs and otherwise we stay at the current one.

Now a brief introduction to the already mentioned Metropolis algorithm is given, 
for details see \cite{rosenthal} or \cite{novak}.
Let $\Omega\subset\R^d$ be a convex body and let $\M=(\Omega,\Le(\Omega),\set{Q(x,\cdot): x\in\Omega})$ 
be a reversible Markov scheme with respect to a distribution $\mu$.
With $\Le(\Omega)$ we denote the Lebesgue $\sigma$-algebra 
of $\Omega$ and $Q(x,A)$ is the transition kernel. 
The aim is to simulate a distribution $\mu_\rho$ on the 
measurable space $(\Omega,\Le(\Omega))$, 
which is defined by an unnormalised density $\rho$ such that
\begin{equation} \label{mur}
\mu_\rho(A)=\frac{\int_A \rho(x)\,\mu(dx) }{\int_\Omega \rho(x)\,\mu(dx)}.
\end{equation}  
It is required that we have an oracle for the evaluation of $\rho$.
In this setup a Metropolis step works as described in Algorithm~\ref{metro_step}.
\noindent
\begin{alg}[ht]
\begin{flushleft}
\textbf{\quad $X_{n+1}=\mbox{metro\_step}(X_n,\rho(\cdot),\,Q(X_n,\cdot))$} 
\end{flushleft}
\hrule
\begin{enumerate}
  \item 	\texttt{choose $Y$ from $Q(X_n,\cdot)$;}\\
  \item   \texttt{calculate $
                            \gamma:=\rho(Y)/\rho(X_n);
                            $}
  \item   \texttt{if $\gamma\geq rand()$ then\\ \indent\quad\, return $Y$;}
  \item   \texttt{else return $X_n$.}
\end{enumerate}
\caption{Metropolis step from $X_n$ to $X_{n+1}$}
\label{metro_step}
\end{alg}
The procedure $rand()$ returns a uniformly distributed random number between zero and one.
If we choose a starting point $X_0$ from a known distribution and take this as input in the method, then we obtain, after repeating Algorithm~\ref{metro_step}, a Markov chain on $\Omega$. The corresponding Markov kernel is defined by
\begin{equation} \label{K_rho}
K_\rho(x,A):=\int_A \theta(x,y)\, Q(x,dy)+ I(x,A)\left( 1-\int_\Omega \theta(x,y)\,Q(x,dy) \right),
\end{equation}
where 
\[
I(x,A)=\begin{cases} 1 & x\in A\\
                     0 & x\in A^c
       \end{cases}
\quad\mbox{and}\quad 
\theta(x,y):=\min\left\{ 1,\frac{\rho(y)}{\rho(x)} \right\}.
\]
The next implication confirms that the resulting Markov scheme 
$\mathcal{M}_\rho= (\Omega,\mathcal{L}(\Omega),\linebreak\{ K_\rho (x,\cdot) : x\in\Omega\})$ 
is reversible concerning $\mu_\rho$. 

\begin{lemma}
If the proposal Markov scheme $\M$ of the Metropolis Hastings method is reversible 
with respect to a distribution $\mu$, 
then the reversibility condition holds also for $\M_\rho$ with respect to $\mur$.
\end{lemma}

\begin{proof}
It is enough to show that the identity 
\[
\int_A K_\rho(x,B)\,\mu_\rho(dx)=\int_B K_\rho(x,A)\,\mu_\rho(dx)
\]
for disjoint sets $A,B\in\Le(\Omega)$ is true.
Furthermore $\theta(y,x)\rho(y)=\theta(x,y)\rho(x)$ for $x,y \in \Omega$
and we define $k:=\int_\Omega \rho(x) \;\mu(dx)$.
Hence this implies
\begin{align*}
&\quad\int_A K_\rho(x,B)\,\mu_\rho(dx) 
\underset{\eqref{K_rho}}{=} \int_A \int_B \theta(x,y)\; Q(x,dy)\,\mu_\rho(dx)\\
&\underset{\eqref{mur}}{=} \frac{1}{k} \int_A \int_B \theta(x,y)\rho(x)\; Q(x,dy)\,\mu(dx)\\
&= \frac{1}{k} \int_\Omega \int_\Omega \chi_A(x)\chi_B(y)\, \theta(x,y)\rho(x)\; Q(x,dy)\,\mu(dx)\\
&\underset{\eqref{F_x_y}}{=} \frac{1}{k}\int_\Omega \int_\Omega \chi_A(y)\chi_B(x)\, \theta(y,x)\rho(y)\; Q(x,dy)\,\mu(dx)\\
&= \frac{1}{k}\int_B \int_A \theta(x,y)\rho(x)\; Q(x,dy)\,\mu(dx)
= \int_B K_\rho(x,A)\,\mu_\rho(dx).
\end{align*}
\end{proof}

Summarizing we have until now a reversible Markov chain on the state space $\Omega$. 
To apply the theory as developed in Section~\ref{error_bounds} the laziness property must be
fulfilled. But as already mentioned we just have to flip a coin and stay at the current state with probability $1/2$. Otherwise do one step with the chain. Formalized written down we consider $\overline{\M}_\rho=(\Omega,\Le(\Omega),\set{\overline{K}_\rho(x,\cdot):x\in\Omega})$, where
\[
\overline{K}_\rho(x,A):=\frac{1}{2}K_{\rho}(x,A)+\frac{1}{2}I(x,A).
\] 
This Markov scheme is lazy, reversible and if it is possible 
to get a lower bound of the conductance we can apply Theorem~\ref{error_thm}. 
Therefore the following result is helpful.
\begin{lemma}  \label{lem_allg_con}
Let $\M=(\Omega,\A,\set{K(x,A):x\in\Omega})$ be an arbitrary reversible Markov scheme concerning $\pi$. 
The conductance of $\overline{\M}=(\Omega,\A,\set{\overline{K}(x,A):x\in\Omega})$, where $\overline{K}(x,A)=\frac{1}{2}K(x,A)+\frac{1}{2}I(x,A)$ is bounded from below, i.e.
\begin{equation*}  \label{allg_conduct}
\phi(\overline{K},\pi)\geq \frac{1}{2} \phi(K,\pi).
\end{equation*}
\end{lemma}
\begin{proof}
The result is obvious after taking the definition of the conductance into account. 
\end{proof}
\begin{rem}
We turned from the Metropolis chain to the lazy one.
Another way would be to ``lazify" the proposal chain 
and after that turn to the Metropolis one.
This is equivalent since
\begin{align*}
K_\rho^{\overline{Q}}(x,A) & = 
\int_A \theta(x,y) \left( \frac{1}{2} Q(x,dy)+\frac{1}{2} I(x,dy) \right)\\
& \quad+I(x,A)\left( 1-\int_\Omega \theta(x,y) \left( \frac{1}{2} Q(x,dy)+\frac{1}{2} I(x,dy) \right)  \right)\\
&=\frac{1}{2}\int_A\theta(x,y)Q(x,dy)+\frac{1}{2}I(x,A)\left( 1-\int_\Omega \theta(x,y) Q(x,dy) \right)+\frac{1}{2}I(x,A)\\
&=\frac{1}{2}K_\rho^Q(x,A)+\frac{1}{2}I(x,A).
\end{align*}
\end{rem}

\subsection{Metropolis algorithm based on the ball walk}
\label{metro_ball_walk}

We come to a concrete given proposal Markov chain, which is defined by a $\d$ ball walk 
on the convex body $\Omega$. This random walk is the same like the already studied one in \cite{novak} and in
different references of volume computation see e.g. \cite{lova_simo1,vempala,vempala_lesson}.
The corresponding Markov scheme is $\M_\d=(\Omega,\Le(\Omega),\set{Q_\d(x,\cdot):x\in\Omega})$,
where
\[
Q_\d(x,A):=\frac{\vol(B(x,\d)\cap A)}{\vol(\d \ball)}+\left( 1-\frac{\vol(B(x,\d)\cap \Omega)}{\vol(\d \ball)} \right)I(x,A).
\]
There $B(x,\d)$ denotes the ball of radius $\d$ around $x\in\Omega$ and $\d\ball:=B(0,\d)$.
We choose $\d\leq D$, where $D$ is the diameter of $\Omega$. 
It is easily seen that $\M_\d$ is reversible concerning the uniform distribution on $\Omega$.
By taking this ball walk as proposal kernel for the Metropolis algorithm 
we get $\M_{\rho,\d}=(\Omega,\Le(\Omega),\set{K_{\rho,\d}(x,\cdot):x\in\Omega})$, where
\[
K_{\rho,\d}(x,A):=\int_A \theta(x,y)\, Q_\d(x,dy)+ I(x,A)\left( 1-\int_\Omega \theta(x,y)\,Q_\d(x,dy) \right).
\]
In \cite{novak} the authors showed that the conductance of the resulting chain is positive
if the density is logconcave and log-Lipschitz. Therefore we consider 
\[
\rad(\Omega):=
\{\rho :\; \rho>0,\; \log\rho\;\; \text{concave},\;  
\abs{\,\log{\rho(x)}-\log{\rho(y)}}\leq \alpha \norm{x-y}{2} \}.
\]
Some more general distributions are studied in \cite{randall_decomp} and \cite{small_world}.
Moreover, let $\Omega$ be the $d$-dimensional unit ball notated by $\ball$ a handy lower bound of the conductance exists.
Thus we can use 
\begin{lemma}  \label{lem_metro_con}
Let the Markov scheme $\M_{\rho,\d}=(\ball,\Le(\ball),\set{K_{\rho,\d}(x,\cdot):x\in\ball})$ be
the Metropolis chain based on the local ball walk $\M_\d$, 
where $\rho\in\rad(\ball)$. 
Then we obtain for an adapted $\d=\min\set{1/\sqrt{d+1},1/\a}$ the following lower bound of the conductance 
\begin{equation} \label{conduct_metro_bsp}
 \phi(K_{\rho,\d},\mur)\geq0.0025\frac{1}{\sqrt{d+1}} \min\set{\frac{1}{\sqrt{d+1}},\frac{1}{\a}}.
\end{equation}  
\end{lemma}
\begin{proof}
See \cite[Corollary~1]{novak}.
\end{proof}
The geometry of the unit ball is essentially used, since the ball walk would get stuck with high probability in domains which have corners. \\


Having finished this we obtain an explicit error bound
of the Markov chain Monte Carlo method on $\Omega:=\ball$ for a
function class $\fad(\ball)$.  
This class is defined by
\begin{equation*} \label{fad}
\fad(\Omega):=\set{(f,\rho):\rho\in\rad(\Omega), \norm{f}{\infty}\leq1}.
\end{equation*}

The method, based on a certain $\d$ ball walk after a burn-in time $n_0$, 
is presented in Algorithm~\ref{mcmc}, where 
$S_{n,n_0}^\d(f,\rho)=\frac{1}{n}\sum_{j=1}^n f(X_{j+n_0})$ if $(f,\rho)\in\fad(\ball)$. 
At first we should care about the starting point in $\ball$.
The simplest way to handle this is choosing the initial state concerning the uniform distribution on the state space $\ball$. 
So the following  calculation for $\nu$, where $A\in\Le(\ball)$  holds true
\[
\nu(A)=\frac{\vol(A)}{\vol(\ball)}=\frac{1}{\vol(\ball)}\int_A \int_\ball \frac{\rho(y)}{\rho(x)} \,dy\, \mur(dx).
\] 
This implies that for $\rho\in\rad(\ball)$
\[
\norm{\frac{d\nu}{d\mur}}{\infty}\leq \exp(2\a).
\]


Now let us turn our view to the error of this Markov Chain Monte Carlo method
and summarize the previous outcomes.

\begin{theorem}  \label{main_appl}
Let $X_1,X_2,\dots$ be the lazy Metropolis Markov chain which is based on a $\d$ ball walk, where $\d=\min\set{1/\sqrt{d+1},1/\a}$. Furthermore it is required that $(f,\rho)\in\fad(\ball)$.
Then we get 
\[
e(S_{n,n_0}^\d,f)\leq 
8000 \frac{\sqrt{d+1} \max\set{ \sqrt{d+1},\a}}{\sqrt{n}},
\]
where $n_0\geq1280000\cdot\a(d+1)\max\set{d+1,\a^2}$.
\end{theorem}
\begin{proof}
After the consideration for the initial distribution $\nu$, the lower bound \eqref{conduct_metro_bsp} for the conductance and  applying Lemma~\ref{lem_allg_con}, Lemma~\ref{lem_metro_con} and \eqref{error_case} the claim is proven.
\end{proof}

For an interpretation let us consider the cost of the underlying method.
With Theorem~\ref{main_appl} we have 
\begin{align*}
 \cost(f,\e)
\leq& \left\lceil 1280000\cdot\a (d+1) \max\set{d+1,\a^2} \right\rceil \\
\notag&+     \left\lceil 64000000\cdot(d+1)\max\set{d+1,\a^2}\e^{-2}\right\rceil.
\end{align*}
This shows that the cost depends only polynomial on the dimension and the Lip\-schitz constant 
such that the suggested algorithm $S_{n,n_0}$ avoids the curse of dimension.
In this setting it is worth to mention that the number of time steps $n+n_0$ is proportional 
to the number of function evaluations of $f$ and $\rho$. We need at most $n+n_0$ oracle calls for $\rho$ and $n$ for $f$.



\subsection*{Acknowledgements}
The author wishes to express his thanks to Erich Novak for 
many suggestions
and several helpful comments 
concerning the presentation.
The author also thanks two anonymous referees for their valuable comments.


\providecommand{\bysame}{\leavevmode\hbox to3em{\hrulefill}\thinspace}
\providecommand{\MR}{\relax\ifhmode\unskip\space\fi MR }
\providecommand{\MRhref}[2]{%
  \href{http://www.ams.org/mathscinet-getitem?mr=#1}{#2}
}
\providecommand{\href}[2]{#2}

\end{document}